\RequirePackage{ifpdf}
\ifpdf 
\documentclass[pdftex]{sigma}
\else
\documentclass{sigma}
\fi

\usepackage{xspace,a4wide}
\usepackage[all]{xy}

\numberwithin{theorem}{section}
\numberwithin{proposition}{section} \numberwithin{lemma}{section}
\numberwithin{corollary}{section}
\numberwithin{definition}{section} \numberwithin{example}{section}
\numberwithin{remark}{section}
\newcommand{\C}{\ensuremath{\mathbb{C}}\xspace}

\renewcommand{\P}{\ensuremath{\mathcal{P}}\xspace}

\newcommand{\al}{\ensuremath{\alpha}}
\renewcommand{\l}{\ensuremath{\lambda}}

\newcommand{\g}{\ensuremath{\mathfrak{g}}}
\newcommand{\G}{\ensuremath{\mathfrak{G}}}
\renewcommand{\H}{\ensuremath{\mathfrak{H}}}
\newcommand{\N}{\ensuremath{\mathfrak{N}}}
\newcommand{\m}{\ensuremath{\mathfrak{m}}}

\newcommand{\B}{\ensuremath{\mathfrak{B}}}

\newcommand{\Z}{\ensuremath{\mathbb{Z}}\xspace}

\newcommand{\R}{\ensuremath{\mathbb{R}}\xspace}

\newcommand{\im}{\operatorname{Im}}

\newcommand{\Span}{\operatorname{Span}\xspace}

\renewcommand{\phi}{\varphi}

\def\D{\Delta}

\begin{document}

\allowdisplaybreaks

\renewcommand{\thefootnote}{$\star$}

\renewcommand{\PaperNumber}{026}

\FirstPageHeading

\ShortArticleName{Induced Modules for Af\/f\/ine Lie Algebras}

\ArticleName{Induced Modules for Af\/f\/ine Lie Algebras\footnote{This paper is a
contribution to the Special Issue on Kac--Moody Algebras and Applications. The
full collection is available at
\href{http://www.emis.de/journals/SIGMA/Kac-Moody_algebras.html}{http://www.emis.de/journals/SIGMA/Kac-Moody{\_}algebras.html}}}

\Author{Vyacheslav FUTORNY and Iryna KASHUBA}

\AuthorNameForHeading{V.~Futorny and I.~Kashuba}

\Address{Institute of Mathematics, University of S\~ao Paulo,\\
Caixa Postal 66281 CEP 05314-970, S\~ao Paulo, Brazil}
\Email{\href{mailto:futorny@ime.usp.br}{futorny@ime.usp.br}, \href{mailto:kashuba@ime.usp.br}{kashuba@ime.usp.br}}

\ArticleDates{Received October 20, 2008, in f\/inal form March 01,
2009; Published online March 04, 2009}

\Abstract{We study  induced modules of nonzero central charge with
arbitrary multiplicities over af\/f\/ine Lie algebras. For a given
pseudo parabolic subalgebra~${\mathcal P}$ of an af\/f\/ine Lie algeb\-ra~${\mathfrak G}$,
 our main result establishes the equivalence between a certain category of
  ${\mathcal P}$-induced ${\mathfrak G}$-modules and the category of weight ${\mathcal P}$-modules with injective action
of the central element of ${\mathfrak G}$. In particular, the induction
functor preserves irreducible
 modules. If ${\mathcal P}$ is a~parabolic subalgebra with a
 f\/inite-dimensional Levi factor then it def\/ines a unique
  pseudo parabolic subalgebra ${\mathcal P}^{ps}$, ${\mathcal P}\subset {\mathcal P}^{ps}$. The
  structure of ${\mathcal P}$-induced modules in this case is fully determined by the
  structure of ${\mathcal P}^{ps}$-induced modules.
 These results generalize similar reductions in particular cases previously considered
 by  V.~Futorny, S.~K\"onig,  V.~Mazorchuk [{\it Forum Math.} {\bf 13} (2001), 641--661],
 B.~Cox [{\it Pacific J. Math.} {\bf 165} (1994), 269--294] and  I.~Dimitrov, V.~Futorny, I.~Penkov [{\it Comm. Math. Phys.} {\bf 250} (2004), 47--63].}


\Classification{17B65; 17B67}

\section{Introduction}\label{s1}

It is dif\/f\/icult to over-estimate the importance of Kac--Moody
algebras for modern mathematics and physics. These algebras were
introduced in 1968 by V.~Kac and R.~Moody as a~generalization of
simple f\/inite-dimensional  Lie algebras, by relaxing the condition
of Cartan matrix to be positive def\/inite. We address to~\cite{K}
for the basics of the Kac--Moody theory.

{\em Affine} Lie algebras are the most studied among
inf\/inite-dimensional Kac--Moody algebras, and have very wide
applications. They correspond to the case of positive semidef\/inite
matrix ($\det (a_{ij})=0$, with positive principal minors).

 All results in
the paper hold for both {\em untwisted} and {\em twisted} af\/f\/ine
Lie algebras of rank greater than $1$. Let $\G$ be an af\/f\/ine
Kac--Moody algebra with
 a $1$-dimensional center
$Z=\C c$.

A natural way to construct representations of af\/f\/ine Lie algebras
is via induction from parabolic subalgebras. Induced modules play
an important role in the classif\/ication problem of irreducible
modules. For example, in the f\/inite-dimensional setting any
irreducible weight module is a quotient of the module induced from
an irreducible module over a~parabolic subalgebra, and this module
is \emph{dense} (that is, it has the largest possible set of
weights) as a~module over the Levi subalgebra of the parabolic
\cite{Fe, F1,DMP}. In particular, dense irreducible
module is always torsion free if all weight spaces are
f\/inite-dimensional. In the af\/f\/ine case, a similar conjecture
\cite[Conjecture 8.1]{F2} singles out induced modules as
construction devices for irreducible weight modules. This
conjecture has been shown for $A_1^{(1)}$ \cite[Proposition~6.3]{F3}  and for all af\/f\/ine Lie algebras in the case of modules with
f\/inite-dimensional weight spaces and nonzero central charge~\cite{FT}. In the latter case, a phenomenon of reduction to
modules over a proper subalgebra (f\/inite-dimensional reductive)
provides a classif\/ication of irreducible modules. Moreover, recent
results of I.~Dimitrov and D.~Grantcharov~\cite{DG} show the
validity of the conjecture also for modules with
f\/inite-dimensional weight spaces and zero central charge.

For  highest weight modules (with respect to nonstandard Borel
subalgebras) with nonzero central charge such reduction was shown
in \cite{C, FS}. The case of induced modules from a~parabolic subalgebra with a f\/inite-dimensional Levi factor was
considered in~\cite{FKM1}. In particular, it was shown that such
categories of modules are related to projectively stratif\/ied
algebras~\cite{FKM2}. A more general setting of toroidal Lie
algebras was considered in~\cite{DFP} for induced modules from
general Borel subalgebras.

The main purpose of the present paper is to show, that in the
af\/f\/ine setting all known cases of the reduction are just
particular cases of a general reduction phenomenon for modules
with nonzero central charge.

 We assign to each
\emph{parabolic subset} $P$ of the root system $\Delta$ of $\G$,
the {\em parabolic subalgebra} $\G_{P}=\G_{P}^0\oplus \G_{P}^+$ of
$\G$ with the Levi subalgebra $\G_{P}^0$ (cf.\ Section~\ref{section3.1}). The
subalgebra $\G_{P}^0$ is inf\/inite-dimensional if and only if
$P\cap -P$ contains  imaginary roots. If, in the same time, $P\cap
-P$ contains some real roots then we def\/ine our key subalgebra
$\G_P^{ps}\subset\G_{P}$, which will be called {\em pseudo
parabolic subalgebra}.

Denote by $U(\G)$  the universal enveloping algebra of $\G$.

Let $\P=\P_0\oplus \N$ be a (pseudo) parabolic subalgebra of $\G$
with the Levi subalgebra $\P_0$. If $N$ is a weight $\P_0$-module
then it can be viewed as a $\P$-module with a trivial action of
$\N$. Then one can construct the induced $\G$-module ${\rm
ind}(\P, \G; N)=U(\G)\otimes_{\P} N$. Hence
\[
{\rm ind}(\P, \G): \ N\longmapsto {\rm ind}(\P, \G; N)
\]
def\/ines a functor from the category of weight $\P_0$-modules to the category of weight $\G$-modules.

Denote by $W(\P_0)$  the full subcategory of weight $\P_0$-modules
on which the central element $c$ acts injectively. Our main result
is the following {\em reduction} theorem.

\begin{theorem}\label{induced-irred}
Let $\G$ be affine Lie algebra of rank greater than $1$,
$\P=\P_0\oplus \N$ a pseudo parabolic subalgebra of $\G$,
$\P_0$ infinite-dimensional and ${\rm ind}_0(\P, \G)$ the
restriction of the induction functor ${\rm ind}(\P, \G)$ onto
$W(\P_0)$. Then the functor ${\rm ind}(\P, \G)$ preserves the
irreducibles.
\end{theorem}

This result allows to construct new irreducible modules for af\/f\/ine
algebras using  parabolic induction from af\/f\/ine subalgebras.

Theorem~\ref{induced-irred} follows from a more  general result (see Theorems~\ref{induced-verma-irred}
and~\ref{theorem-equivalence}) which establishes an equivalence of certain categories of
modules.

In the case when the Levi factor $\G_{P}^0$ is f\/inite-dimensional,
we def\/ine a certain  subalgebra $\m_{P} \subset \G$, which leads
to a parabolic decomposition $\G = \N_{P}^- \oplus \m_{P} \oplus
\N_{P}^+$ with $\N_{P}^+\subset \G_{P}^+$ and $\G_{P}^0\subset
\m_{P}$ (cf.\ Section~\ref{section2}).
 If $\m_{P}$ is f\/inite-dimensional then
$\G_{P}^0=\m_{P}$ and $\G_{P}^+=\N_{P}^+$. In this case there is
no reduction. If $\m_{P}$ is inf\/inite-dimensional then $\m_{P}
\oplus \N_{P}^+$ is pseudo parabolic subalgebra and the reduction
theorem above implies Theorem~8 in~\cite{FKM1}.

The structure of the paper is as follows. In Section~\ref{section2} we recall
the classif\/ication of Borel subalgebras and parabolic subalgebras.
Section~\ref{section3} is devoted to the study of parabolic and pseudo
parabolic induction. In particular, we prove Theorem~\ref{induced-verma-irred} that
describes the structure of induced modules. In the last section we
introduce certain categories of $\G$-modules and of $\P$-modules
and establish their equivalence (Theorem~\ref{theorem-equivalence}).

\section{Preliminaries on Borel subalgebras\\ and parabolic subalgebras}\label{section2}

Let $\H$ be a Cartan subalgebra of $\G$ with the following
decomposition
\[
\G = \H \oplus (\oplus_{\alpha \in \H^* \backslash \{0\}}
\G_\alpha),
\]
where $\G_\alpha := \{ x \in \G \, | \, [h, x] = \alpha(h) x
{\text{ for every }} h \in \H\}$. Denote by \mbox{$\Delta = \{ \alpha
\in \H^* \backslash \{0\} \, | \, \G_\alpha \neq 0\}$}  the root
system of $\G$. Let $\pi$ be a basis of the root system $\D$.
 For a subset
$S\subset \pi$  denote by~$\D_+(S)$  the subset of~$\D$,
consisting of all linear combinations of elements of $S$ with
nonnegative coef\/f\/icients. Thus,
 $\D_+(\pi)$ is the set of positive roots
with respect to~$\pi$. Let $\delta\in \D_+(\pi)$ be the
indivisible imaginary root.
  Then the set of all imaginary roots is $\Delta^{im}=\{k\delta|
k\in\Z\setminus\{0\}\}$. Denote  by $G$ a Heisenberg subalgebra of
$\G$ generated by the root spaces $\G_{k\delta}$, $k\in \Z\setminus \{0\}$.

Denote by $S_{\pi}$ a root subsystem generated by $S$ and
$\delta$. Let $S_{\pi}^+=S_{\pi}\cap \D_+(\pi)$. For a subset
$T\subset \D$ denote by $\G(T)$ the subalgebra of $\G$ generated
by the root subspaces $\G_{\alpha}$,
 $\alpha\in T$, and  let $\H(T)=\H\cap \G(T)$. The subalgebra $\G(-T)$ will be called the {\em opposite subalgebra}
 to $\G(T)$.

 Let $S\subset \pi$, $S=\cup_i
S_i$ where all  $S_i$'s are connected and $S_i\cap S_j=\varnothing$ if $i\neq j$.

\begin{proposition}[\protect{\cite[Proposition~2]{FKM1}}]\label{prop-decomp-Levi}
$\G(S_{\pi})=\G^S + G(S) +\H$, where $\G^S=\sum_{i}\G^i$, $[\G^i,
\G^j]=0$, $i\neq j$, $\G^i$ is the derived algebra of an affine
Lie algebra of rank $|S_i|+1$, $[\G^S, G(S)]=0$, $G(S)\subset G$,
$G(S)+(G\cap \G^S)=G$, $\G^S\cap G(S)=\cap_{i}\G^i=Z$.
\end{proposition}

Let $V$ be a weight $\G$-module, that is $V=\oplus_{\mu\in \H^*}
V_{\mu}$, \mbox{$V_{\mu}=\{v\in V\, |\, hv=\mu(h)v, \forall \, h\in \H\}$}. If $V$
is irreducible then $c$ acts as a scalar on $V$, it is called the
{\em central charge} of $V$. A~classif\/ication of irreducible
modules in the category  of all weight modules is still an open
question even in the f\/inite-dimensional case. In the af\/f\/ine case,
a classif\/ication is only known in the subcategory of modules with
f\/inite-dimensional weight spaces (\cite{FT} for nonzero charge and
\cite{DG} for zero charge), and in certain subcategories of
induced modules with some inf\/inite-dimensional weight spaces~\cite{F2, FKM1}. If $V$ is a weight module (with respect to
a f\/ixed Cartan subalgebra) then we denote by $w(V)$ the set of
weight, that is $w(V)=\{\lambda\in \H^*| V_{\lambda}\neq 0\}$.

The af\/f\/ine Lie algebra $\G$
has the associated simple f\/inite-dimensional Lie algebra $\g$ (for
details see \cite{K}). Of course, in the untwisted
case $\G$ is just the af\/f\/inization of $\g$:
\[
\G=\g\otimes \C\big[t,t^{-1}\big]\oplus\C c\oplus \C d,
\]
where $d$ is the degree derivation: $d(x\otimes t^n)=n(x\otimes t^n)$, $d(c)=0$, for all $x\in \G$,
$n\in \Z$.

A closed subset $P\subset \D$ is called a {\em partition} if
$P\cap (-P)=\varnothing$ and $P\cup (-P)=\D$. In the case of
f\/inite-dimensional simple Lie algebras,  every partition
corresponds to a choice of positive roots in $\D$, and all
partitions are conjugate by the Weyl group. The situation is
dif\/ferent in the inf\/inite-dimensional case. In the case of  af\/f\/ine
Lie algebras the partitions are divided into a f\/inite number of
Weyl group orbits (cf.~\cite{JK,F2}).

Given a partition $P$ of $\D$, we def\/ine a {\em Borel} subalgebra
$\B_P\subset \G$ generated by $\H$ and the root spaces ${\mathfrak
G}_{\alpha}$ with $\alpha \in P$. Hence, in the af\/f\/ine case not
all of the Borel subalgebras are conjugate but there exists a
f\/inite number of conjugacy classes.

A {\em parabolic} subalgebra of $\G$  corresponds to a {\em
parabolic} subset $P\subset \D$, which is a closed subset in $\D$
such that $P\cup (-P)=\D$. Given such a parabolic subset $P$, the
corresponding parabolic subalgebra $\G_{P}$ of $\G$ is generated
by $\H$ and all the root spaces ${\mathfrak G}_{\alpha}$, $\alpha
\in P$.

The conjugacy classes of Borel subalgebras of $\G$ are
parameterized by the parabolic subalgebras of the associated
f\/inite-dimensional simple Lie algebra $\g$. We just recall the
construction in the untwisted case. Parabolic subalgebras of $\g$
are def\/ined as above, they correspond to parabolic subsets of the
roots system of $\g$. Let $\mathfrak p=\mathfrak p_0\oplus
\mathfrak p_+$ be a parabolic subalgebra of $\g$ containing a
f\/ixed Borel subalgebra $\mathfrak b$ of $\g$. Def\/ine
\[
\B({\mathfrak p})=\mathfrak p_+\otimes \C[t,t^{-1}]\oplus \mathfrak p_0\otimes t\C[t]
\oplus \mathfrak b \oplus \C c\oplus \C d.
\] For any Borel
subalgebra $\B$ of $\G$ there exists a parabolic
subalgebra $\mathfrak p$ of $\g$ such that $\mathfrak B$ is
conjugate to $\B({\mathfrak p})$~\cite{JK,F2}.

Any Borel subalgebra conjugated to  $\B({\g})$ is called {\em
standard}. It is determined by a choice of positive roots in $\G$.
Another extreme case  $\mathfrak p_0=\H$ corresponds to the {\em
natural} Borel subalgebra  of $\G$.

We will also  use the  geometric description  of Borel subalgebras
in $\G$ following \cite{DP}. Let $W=\Span_\R{\Delta}$, $n=\dim
W$. Let
\[
F = \{ \{0\} = F_n \subset F_{n-1} \subset
\cdots \subset F_1 \subset F_0 = W \}
\] be a f\/lag of maximal
length in~$W$. The f\/lag $F$ of maximal length is called oriented
if, for each~$i$, one of the connected components of $F_i
\backslash F_{i+1}$ is labelled by $+$ and the other one is
labelled by~$-$. Such an oriented f\/lag $F$ determines the
partition
\[
P:= \D\cap \big(\cup_i (F_i \backslash F_{i+1})^+\big)
\]
 of $\D$. Denote
$P_i^0=\D\cap F_i$. The subsets $P_i^0$ are important invariants
of the partition $P$.

Next statement follows immediately from the description of partitions of root systems~\cite{F4}.

\begin{proposition} \label{prop1}
Let $P\subset \D$ be  a partition.
 There exists an oriented flag of maximal length $F =
\{ \{0\} = F_n \subset F_{n-1} \subset \cdots \subset F_1 \subset
F_0 = W \}$ which determines the partition $P$.
\end{proposition}

Given a partition $P$ of $\D$ we will denote by $F(P)$ the
corresponding oriented f\/lag of maximal length.

Let $\lambda: B_{P}\rightarrow \C$ be a $1$-dimensional
representation of $B_{P}$. Then one def\/ines an induced {\em Verma
type} ${\G}$-module
\[
M_{P}(\lambda)=U({\G})\otimes_{U(B_{
P})}\C.
\] The module $M_{\Delta_+(\pi)}(\lambda)$ is a classical Verma
module with the highest weight $\lambda$~\cite{K}. In the case of
natural Borel subalgebra we obtain {\em imaginary} Verma modules
studied in~\cite{F5}.

Let $\B_P=\H\oplus \N_P$, where $\N_P$ is generated by ${\mathfrak
G}_{\alpha}$, $\alpha \in P$. Note that the module $M_{
P}(\lambda)$ is $U(\N_P^{-})$-free, where $\N_P^{-}$ is the
opposite subalgebra to $\N_P$. The theory of Verma type modules
was developed in~\cite{F2}. It follows immediately from the
def\/inition that Verma type module with highest weight $\lambda$
has a unique maximal submodule. Also, unless it is a classical
Verma module, it has both f\/inite and inf\/inite-dimensional weight
spaces and it can be obtained using the parabolic induction from a
classical Verma module $M$ with highest weight $\lambda$ over a
certain inf\/inite-dimensional Lie subalgebra. Moreover, if the
central charge of such Verma type module is nonzero, then the
structure of this module is completely determined by the structure
of module~$M$, which is well-known~\cite{F2, C,FS}.

Let $P\subset \D$ be a parabolic subset, $P\cup (-P)=\D$, $\G_{P}$
 the corresponding parabolic subalgebra  of $\G$. Set $P^0=P\cap
-P$. Then
\[
\G=\G_{P}^-\oplus \G_{P}^0\oplus \G_{P}^+,
\] where
$\G_{P}^{\pm}=\sum_{\al\in P\setminus (-P)}\G_{\pm \al}$ and  $\G_{\P}^0$ is
generated by $\H$ and the subspaces
 $\G_{\al}$ with $\al\in {P}^0$.
 The subalgebra $\G_{P}^0$ is the {\em Levi factor} of $\G$.
A classif\/ication of parabolic subsets of $\D$ and parabolic subalgebras of~$\G$ was obtained in \cite{F2,F4}.

Every parabolic subalgebra $\P\subset \G$ containing a Borel
subalgebra $\B$ has a Levi decomposition $\P=\P_0\oplus \N$, with
the Levi factor $\P_0$ and $\N\subset \B$. Following \cite{F2} we
say that $\P$ has \emph{type I} if $\P_0$ is f\/inite-dimensional
reductive Lie algebra, and $\P$ has {\em type~II} if $\P_0$
contains the Heisenberg subalgebra $G$, generated by the imaginary
root spaces of $\G$.  In the latter case, $\P_0$~is
 an extension of a sum of some af\/f\/ine Lie subalgebras
by a central subalgebra and by a~certain subalgebra of $G$~\cite{F1}.  Note that the radical~$\N$ is solvable only for the
type~II parabolic subalgebras. Type~I parabolic subalgebras are
divided also into two essentially dif\/ferent types depending on
whether $\N$ belongs to some standard Borel subalgebra ({\em type
Ia}) or not ({\em type~Ib}).

It is easy to see that there are parabolic subalgebras $\P$ which
do not correspond to any triangular decomposition of~$\G$~\cite{F2}. In fact, if $\P$ contains $\G_{\al+k\delta}$ for some
$\al$ and inf\/initely many both positive and negative integers $k$
then this parabolic subalgebra does not correspond to any
triangular decomposition of~$\G$. In particular, this is always the case for parabolic subalgebras of type II.

If $\P=\P_0\oplus \N$ is a parabolic subalgebra
of type II then as soon as $\N$ contains $\G_{\al+k\delta}$ for some real root
$\al$ and some $k\in \Z$, it also contains $\G_{\al+r\delta}$ for all $r\in \Z$.

Geometrically parabolic subsets correspond to partial oriented
f\/lags of maximal length. Let $P$ be a parabolic subset which
contains a partition $\tilde{P}$ and $F(\tilde{P})$ the
corresponding full oriented f\/lag: $F(\tilde{P})=\{ \{0\} = F_n
\subset F_{n-1} \subset \cdots \subset F_1 \subset F_0 \}.$ If
$P^0\subset F_k$ for some $k$ and $P^0$ is not in $F_{k+1}$ then
$P$ is completely determined by the partial oriented f\/lag
\[
F(P)=\{ \{0\} = F_k\subset
\cdots \subset F_1 \subset F_0=W \}.
\]
Here $P^0=F_k\cap \Delta$, $\P_0=\G_P^0$.
The corresponding parabolic subalgebra $\P$  has type II if $\Delta^{\im}\subset F_k$ and it has type~Ib if
$\Delta^{\im}\subset F_s$ for some $s$, $1\leq s< k$.

Let $\P=\P_0\oplus \N$ be a non-solvable parabolic subalgebra  of
type II, $\P_0=[\P_0, \P_0]\oplus G(\P)+\H$, where $G(\P)\subset
G$ is the orthogonal completion (with respect to the Killing form)
of the Heisenberg subalgebra  of $[\P_0, \P_0]$, that is $G(\P)+
([\P_0, \P_0]\cap G)=G$ and $G(\P)\cap [\P_0, \P_0]=\C c$. Note
that by Proposition~\ref{prop-decomp-Levi}, $[\P_0, \P_0]$ is a
sum of af\/f\/ine subalgebras of $\G$. Let
\[
G(\P)=G(\P)_-\oplus \C
c\oplus G(\P)_+
\] be a triangular decomposition of $G(\P)$. Def\/ine
a {\rm pseudo parabolic subalgebra} $\P^{ps}=\P^{ps}_0\oplus
\N^{ps}$, where $\P^{ps}_0$ is generated by the root spaces
$\G_{\alpha}$, $\alpha\in P\cap -P\cap \D^{\rm re}$ and $\H$,
while  $\N^{ps}$ is generated by the root spaces $\G_{\alpha}$,
$\alpha\in P\setminus (-P)$ and $G(\P)_+$.  Then $\P^{ps}$ is a
proper subalgebra of $\P$.

Suppose $\P=\P_0\oplus \N$ is a parabolic subalgebra
of type Ib, $P$ the corresponding parabolic subset and $F(P)$ the partial f\/lag. Then $\P_0$
is a f\/inite-dimensional reductive Lie
algebra.
Assume
 $\G_{\delta}\subset \N$ and $\Delta^{\im}\subset F_s$, with the largest such $s$, $1\leq s< k$. Note that
 for any  $\al\in P^0$,  $P\setminus P^0$ contains the roots of the form
$\al+k\delta$ and $-\al+k\delta$ for all $k>0$.
 Denote by $\m_{P}$ a~subalgebra of $\G$ generated by  $F_s\cap \D$ and $H$. This is an inf\/inite-dimensional
 Lie algebra which contains $\P_0$ and  $\tilde{\m}_P=\m_P\cap \P$ is a parabolic subalgebra of $\m_P$ of type Ia.

 Let
 $N_P$ be the span of all root subspaces $\G_{\beta}$, $\beta\in P$ which are not in $\m_P$, $N_P\subset \N$. Then
   $\P=\tilde{\m}_P\oplus N_P$. It follows immediately that

\begin{proposition}[\protect{\cite[Proposition~3]{FKM1}}]\label{pp2}
 $\tilde{\m}_P=\m_P \oplus
N_P$ is a parabolic subalgebra of $\G$ of type II and $\G=N_P^{-}\oplus \m_P \oplus
N_P$, where $N_P^{-}$ is the opposite algebra to $N_P$.
\end{proposition}

Hence, any parabolic subalgebra of type Ib can be extended
canonically to the parabolic subalgebra of type II. Moreover, it
can be extended canonically to the pseudo parabolic subalgebra
$\tilde{\m}^{ps}_P$:
\[
\P\subset \tilde{\m}^{ps}_P\subset \tilde{\m}_P.
\]

\section{Parabolic induction}\label{section3}

Let $P$ be a parabolic subset of $\D$.
 Let $N$ be a weight (with respect to $\H$) module over the
 parabolic subalgebra
 $\P=\G_{P}$ (respectively pseudo parabolic subalgebra $\P^{ps}$),    with a
trivial action of $\G_{P}^+$  (respectively  $(\G_{P}^{ps})^+$),
and let
\[
M_{P}(N)={\rm ind}(\G_{P}, \G; N), \qquad M_{P}^{ps}(N)={\rm ind}(\G_{P}^{ps}, \G; N)
\]
 be the induced $\G$-modules. If $N$ is irreducible then  $M_{P}(N)$ (respectively $M_{P}^{ps}(N)$)  has a~unique irreducible quotient
$L_{P}(N)$ (respectively $L_{P}^{ps}(N)$). If $N$ is irreducible
$\P$-module such that $G(\P)_+$ acts trivially on $N$ then
$M_{P}(N)\simeq M_{P}^{ps}(N)$ and $L_{P}(N)\simeq L_{P}^{ps}(N)$.

If $\G_{P}^0\neq \G$ then $L_{P}(N)$ is said to be {\em
parabolically induced}. Following~\cite{DMP} we will say that
irreducible $\G$-module $V$ is {\em cuspidal} if it is not of type
$L_{P}(N)$ for any proper parabolic subset $P\subset \D$ and any~$N$.

We see right away that a classif\/ication  of  irreducible   $\G$-modules reduces to the classif\/ication  of all
irreducible cuspidal modules over Levi subalgebras of $\G$. Namely we have

\begin{proposition}
Let $V$ be an irreducible weight $\G$-module. Then there exists
a parabolic subalgebra $\P=\P_0\oplus \N$ of $\G$ $($possibly equal $\G)$ and an irreducible
weight cuspidal $\P_0$-module~$N$ such that
$V\simeq L_P(N)$.
\end{proposition}

\begin{example}\qquad
\begin{itemize}\itemsep=0pt
\item Let $V$ be an irreducible weight cuspidal $\g$-module then
$V\otimes \C[t,t^{-1}]$ is an irreducible cuspidal $\G$-module
with zero central charge.

\item Modules obtained by the parabolic
induction from  cuspidal modules over the Heisenberg subalgebra
are called {\em loop} modules~\cite{Ch}.

\item Pointed (that is,
all weight spaces are  $1$-dimensional) cuspidal modules were
studied in~\cite{S}.
\end{itemize}
\end{example}

A Levi subalgebra of $\G$ is {\em cuspidal} if it admits a weight
cuspidal module. All cuspidal Levi subalgebras of type~Ia and~Ib
parabolics were classif\/ied in \cite{Fe}. They are the subalgebras
with simple components of type~$A$ and~$C$. All cuspidal Levi
factors of  type~II parabolic subalgebras were described in
\cite{F4,F2}. For any af\/f\/ine Lie algebra the simplest Levi
subalgebra of type~II is a Heisenberg subalgebra. Below we provide
a list of all other Levi subalgebras of type~II with the connected
root system.

\begin{center}
\begin{tabular}{|l|l|}
\hline $\G$ & $\P_0$ \rule{0pt}{14pt}\\ \hline
$\,A^{(1)}_n\,$ & $A^{(1)}_k,$ $1\leq k\leq n-1$ \rule{0pt}{14pt}\\
\hline
$\,B^{(1)}_n\,$ & $A^{(1)}_k,$ $1\leq k\leq n-1$, $\ C_2^{(1)}$, $\ B_k^{(1)}$, $3\leq k\leq n-1$ \rule{0pt}{14pt}\\
\hline
$\,C^{(1)}_n\,$ & $A^{(1)}_k,$ $1\leq k\leq n-1$, $\ C_k^{(1)}$, $2\leq k\leq n-1$ \rule{0pt}{14pt}\\
\hline
$\,D^{(1)}_n\,$ & $A^{(1)}_k,$ $1\leq k\leq n-1$, $\ D_k^{(1)}$, $4\leq k\leq n-1$ \rule{0pt}{14pt}\\
\hline
$\,G^{(1)}_2,\,$ $D_4^{(3)}$ & $A^{(1)}_1$ \rule{0pt}{14pt}\\
\hline
$\,F^{(1)}_4\,$ & $A_1^{(1)}$, $A_2^{(1)}$, $C_2^{(1)}$, $C_3^{(1)}$, $B_3^{(1)}$  \rule{0pt}{14pt}\\
\hline $\,E^{(1)}_l,\,$ $l=6,7,8\,$ & $A^{(1)}_k,$ $1\leq k\leq
l-1$, $\ D^{(1)}_k,$ $4\leq k\leq l-1$,
$\ E^{(1)}_k,$ $6\leq k\leq l-1$  \rule{0pt}{14pt}\\
\hline $\,A^{(2)}_{2n}\,$ & $A^{(1)}_k,$ $1\leq k\leq n-1$, $\
A^{(2)}_{2k},$ $1\leq k\leq n-1$,
$\ E^{(1)}_k,$ $6\leq k\leq l-1$ \rule{0pt}{14pt}\\
\hline $\,D^{(2)}_n\,$  & $A^{(1)}_k,$ $1\leq k\leq n-2$, $\
D^{(2)}_{k},$ $3\leq k\leq n-1$,
$\ E^{(1)}_k,$ $6\leq k\leq l-1$ \rule{0pt}{14pt}\\
\hline $\,A^{(2)}_{2n-1}\,$ & $A^{(1)}_k,$ $1\leq k\leq n-2$, $\
A^{(2)}_{2k-1},$ $3\leq k\leq n-1$,
$\ D^{(2)}_3$ \rule{0pt}{14pt}\\
\hline
$\,E^{(2)}_6\,$ & $A_1^{(1)}$, $A_2^{(1)}$, $D_3^{(2)}$, $D_4^{(2)}$, $D_5^{(2)}$\rule{0pt}{14pt}\\
\hline
\end{tabular}
\end{center}

 A nonzero element $v$
of a $\G$-module $V$ is called $\P$-{\it primitive} if
$\G_{P}^+v=0$.
   Let $Q_{P}$ be the free Abelian group
generated by $P^0$. The following statement is standard.

\begin{proposition}\label{pr2} Let  $V$ be an irreducible weight $\G$-module
with a $\P$-primitive element of weight $\l$, $P^0\neq \Delta$,
 $N=\sum_{\nu \in Q_{P}}V_{\l+\nu}$. Then  $N$ is an
irreducible $\G_{P}$-module and $V$ is
isomorphic to $L_{P}(N)$.
\end{proposition}

 If $V$ is generated by a $\B$-primitive
element $v\in V_{\l}$ then $V$ is a highest weight module with
highest weight $\l$. If $\P$ is of type Ia then $M_{P}(N)$ is a
{\em generalized Verma module} \cite[Section~2]{FKM1}. A
classif\/ication of all irreducible $N$  with f\/inite-dimensional
weight spaces (and hence of $L_{P}(N)$ if $\P$ is of type I) is
known due to Proposition~\ref{pr2}, \cite{M} and \cite{Fe}. Also a
classif\/ication is known when $\P$ is of type~II, $N$ has
f\/inite-dimensional weight spaces and a nonzero charge~\cite{FT}.
In this case, $N$ is the irreducible quotient of ${\rm ind}(\P_0,
\P'; N')$, where $\P'=\P_0'\oplus \N'$ is a parabolic subalgebra
of $\P_0$ of type Ia and $\N'$ is an irreducible $\P'$-module with
a trivial action of $\N'$.

\subsection{Reduction theorem for type II}\label{section3.1}

Let $\pi$ be a basis of the root system $\D$, $\alpha_0\in \pi$ such that
$-\alpha+\delta\in \sum_{\beta\in \pi\setminus \{\alpha_0\}}\Z\beta$ and either
$-\alpha+\delta\in \Delta$ or $1/2(-\alpha+\delta)\in \Delta$.
Let
$\dot{\pi}=\pi\setminus \{\alpha_0\}$ and $\dot{\D}_+$ the free semigroup generated by~$\dot{\pi}$.
Choose a proper subset $S\in \dot{\pi}$ and the root subsystem $S_{\pi}$ generated by $S$ and $\delta$. Set
\[
P_+=\{\alpha+n\delta|\alpha\in \dot{\D}_+\setminus S_{\pi}, n\in \Z\}\cap \D.
\] Then
$P_S=S_{\pi}\cup P_+$ is a parabolic subset with $S_{\pi}=P_S\cap -P_S$. Let
$\G_{P_S}$ be the corresponding parabolic subalgebra. Then it is of type~II
with
\[
\G_{P_S}^0=\G(S_{\pi})=\sum_{\alpha\in S_{\pi}}\G_{\alpha}\oplus \H.
\]

\begin{proposition}[\cite{F2}]\label{prop-parabolic-conjugate}
If $\P$ is a parabolic subalgebra of $\G$ of type II then there
exist $\pi$, $\alpha_0$ and $S$ as above such that $\P$ is
conjugate to $\G_{P_S}$.
\end{proposition}

Hence, it suf\/f\/ices to consider the parabolic subalgebras of type
II in the form $\G_{P_S}$.

Let $S_{\pi}=\cup_i S_i$ be the decomposition of $S_{\pi}$ into
connected components and let $\G_{P_S}^0=\sum_i \G_i\oplus G(P_S)$
be the corresponding decomposition of $\G_{P_S}^0$ (see
Proposition~\ref{prop-decomp-Levi}).

\begin{theorem}\label{induced-verma-irred}
Let $\G$ be of rank $> 1$, $P$  a  parabolic subset of $\D$ such
that $P\cap -P$ contains real and imaginary roots simultaneously.
Consider a weight $\G_{P}^{ps}$-module $V$ which is annihilated by~$(\G_{P}^{ps})^+$ and on which the central element $c$ acts
injectively $($object of $W(\G_P^{ps}))$. Then for any submodule $U$
of ${\rm ind}(\G_{P}^{ps}, \G; V)$ there exists a submodule $V_U$
of $V$ such that
\[
U\simeq {\rm ind}(\G_{P}^{ps}, \G; V_U).
\]
In particular, ${\rm ind}(\G_{P}^{ps}, \G; V)$ is irreducible if
and only $V$ is irreducible.
\end{theorem}

\begin{proof}
The proof follows general lines of the proof of Lemma 5.4 in
\cite{F2}. Denote ${M}^{ps}(V)={\rm ind}(\G_{P}^{ps}, \G; V)$ and
$\hat{M}^{ps}(V)=\sum_{\nu \in Q_{P}, \lambda\in
w(V)}M^{ps}(V)_{\l+\nu}=1\otimes V.$ Then $\hat{M}^{ps}(V)$ is a
$\G_{P}$-submodule of $M^{ps}(V)$ isomorphic to $V$, which
consists of $\P$-primitive elements.

Let $\G_{P}^{ps}=\P_0\oplus \N$, $\N_-$ is the opposite subalgebra
to $\N$.

 Let $U$ be a nonzero submodule of ${M}^{ps}(V)$ and $v\in U$ a
nonzero homogeneous element. Then
\[
v=\sum_{i\in I} u_iv_i,
\]
where $u_i\in U(\N_-)$ are linearly independent homogeneous,
$v_i\in M^{ps}(V)$.

Given a root $\phi\in \D$ denote by ${\rm ht}(\phi)$ the number of
simple roots of $\N_-$ in the decomposition of $\phi$ and by ${\rm
ht_1}(\phi)$ the number of all simple roots  in the decomposition
of $\phi$. Suppose $u_i\in U(\N_-)_{-\phi_i}$. We can assume that
all $\phi_i$'s have the same ${\rm ht}$. Let ${\rm ht}(\phi_i)=1$
for all $i$. Choose $i_0$ such that ${\rm ht_1}(\phi_{i_0})$ is
the least possible. Then there exists a nonzero $x\in \N$ such
that $0\neq xv\in M^{ps}(V)$ and $[x,u_{i_0}]\in U(G\cap \N_-)$
since $[\N, \N_-]\cap G=G(\G_{P}^{ps})$. But $U(G\cap \N_-)$ is
irreducible $G(\G_{P}^{ps})$-module. Hence, there exists $y\in
U(\N)$ such that $y[x, u_{i_0}]v_{i_0}=v_{i_0}$. In the same time
$yxu_iv_i=0$ if $i\neq i_0$. Thus, we obtain $v_{i_0}\in U$ and
$v-u_{i_0}v_{i_0}\in U$. Applying the induction on $|I|$ we conclude
that $v_i\in U$ for all $i\in I$. This completes the proof in the
case ${\rm ht}(\phi)=1$. The induction step is considered
similarly. Hence, $U$ is generated by $U\cap {M}^{ps}(V)$ which
implies the statements.
\end{proof}

\begin{corollary}\label{cor111}
For each $i$, let $V_i$ be an irreducible $G_i$-module with a
nonzero action of the central element. Then ${\rm ind}(\G_P^{ps},
\G; \otimes_i V_i)$ is irreducible.
\end{corollary}

\begin{proof}
Since $\otimes_i V_i$ is irreducible $\G_{P}^{ps}$-module, the
statement follows immediately from
Theorem~\ref{induced-verma-irred}.
\end{proof}

\begin{corollary}\label{corol-nonzero-GVM}
Let $V$ be an irreducible weight non-cuspidal $\G$-module with an
injective action of the central element $c$. Then there exists a
parabolic subalgebra $\P=\P_0\oplus \N$ of $\G$ and an irreducible
weight cuspidal $\P_0$-module $N$ such that $V\simeq L_P(N)$,
where $P$ is the corresponding parabolic subset of $\D$. Moreover,
$V\simeq M_P(N)$ if $\P_0$ is infinite-dimensional and $\N^{ps}
N=0$.
\end{corollary}

\begin{proof}
First statement is obvious. If $\N^{ps} N=0$ then $M_P(N)\simeq
M_P^{ps}(N)$ which is irreducible by
Theorem~\ref{induced-verma-irred}.
\end{proof}

If $\P_0=\H$ then Corollary~\ref{corol-nonzero-GVM} implies reduction theorem
for Verma type modules \cite{C, FS}.

Note that in general ${\rm ind}(\G_P, \G; N)$ need not be
irreducible if $N$ is irreducible. On the other hand we have

\begin{corollary}
For each $i$, let $V_i$ be an irreducible $G_i$-module with the
action of the central element by a nonzero scalar $a$,  $V$ an
irreducible highest weight $G(\G_{P})$-module with highest weight
$a$. Then
\[
{\rm ind}(\G_P, \G; \otimes_i V_i\otimes V)
\] is
irreducible.
\end{corollary}

\begin{proof}
Note that $V$ is isomorphic to the Verma module with highest weight
$a$. Then
\[
{\rm ind}(\G_P, \G; \otimes_i V_i\otimes V)\simeq {\rm ind}(\G_P^{ps}, \G; \otimes_i
 V_i)
 \] which is irreducible by Corollary~\ref{cor111}.
\end{proof}

\begin{corollary} Let $\lambda\in \H^*$, $\lambda(c)\neq 0$, $\B$ a non-standard Borel subalgebra of $\G$, $P$
corresponding partition of $\D$ and
\[
F=\{ \{0\} = F_n\subset
\cdots \subset F_1 \subset F_0=W \}
\] the corresponding full flag of maximal length.
Suppose
 $\delta\in F_{s-1}\setminus F_{s}$ for some $s$, $1\leq s< n$. Denote $\m_s$ the Lie subalgebra of $\G$ generated by the
 root subspaces with roots in $\Delta\cap F_{s-1}$ and by $\H$.   Then $\m_s$ is infinite-dimensional
 and $M_P(\lambda)^s=U(\m_s)v_{\lambda}$ is a highest weight module over $\m_s$. Moreover,
$M_{P}(\lambda)$ is irreducible if and only if $M_P(\lambda)^s$ is
irreducible.
\end{corollary}

\begin{proof}
Indeed, the f\/lag $F$ def\/ines a parabolic subalgebra $\P$ of $\G$
whose  Levi subalgebra is~$\m_s$. Hence,
\[
M_{P}(\lambda)\simeq
{\rm ind}(\P,\G; M_P(\lambda)^s).
\] The statement follows
immediately from  Corollary~\ref{corol-nonzero-GVM}.
\end{proof}

\subsection {Reduction theorem for type Ib}\label{section3.2}

 Let $\P=\P_0\oplus \N$ be a parabolic subalgebra of $\G$ of type Ib, $\dim \P_0<\infty$, $P$ corresponding parabolic subset of $\D$ and
\[
F=\{ \{0\} = F_k\subset
\cdots \subset F_1 \subset F_0=W \}
\] the corresponding partial
f\/lag of maximal length. Then $\P_0$ has the root system $F_k\cap
\Delta$. Since $\P$ is of type~Ib, there exists $s$, $1\leq s< k$,
such that $\delta\in F_s$ and $\delta\notin F_{s+1}$. Then, a Lie
subalgebra  $\m_{P}$ of $\G$ generated by $\H$ and the root spaces
corresponding to the roots from $F_s\cap \D$, is
inf\/inite-dimensional. Obviously, it can be extended to a parabolic
subalgebra $\m_P\oplus \N_P$ of $\G$ of type II, where $\m_{P}$ is
the Levi subalgebra.

\begin{corollary}\label{cor-induced-irred-type Ib}
Let $P$ be a parabolic subset of $\D$ such that $\G_{P}^0$ is
finite-dimensional and $\m_{P}$ is infinite-dimensional. Consider
a weight $\m_{P}$-module $V$ which is annihilated by $\N_{P}^{ps}$
and on which the central element $c$ acts injectively. Then for
any submodule $U$ of ${\rm ind}(\m_{P}\oplus \N_P, \G; V)$ there
exists a submodule $V_U$ of $V$ such that
\[
U\simeq {\rm ind}(\m_{P}\oplus\N_P, \G; V_U).
\]
In particular, ${\rm ind}(\m_{P}\oplus\N_P, \G; V)$ is irreducible
if and only if $V$ is irreducible.
\end{corollary}

\begin{proof}
Consider the pseudo parabolic subalgebra $\m_{P}^{ps}\oplus
\N_P^{ps}\subset \m_{P}\oplus\N_P$, $\G_{P}\subset
\m_{P}^{ps}\oplus \N_P^{ps}$. Then
\[
 {\rm ind}(\m_{P}\oplus \N_P, \G; U(\m_P)V)\simeq
{\rm ind}(\m_{P}^{ps}\oplus \N_P^{ps}, \G; U(\m_P^{ps})V),
\] since
$U(\m_P^{ps})V\simeq {\rm ind}(\G_P^0, \m_P^{ps}; V)$ and
$U(\m_P)V\simeq U(\m_P^{ps})V\otimes M,$ where $M$ is the hig\-hest~weight $G(\G_P)$-module of highest weight $a$. Hence, the
statement follows from Theorem~\ref{induced-verma-irred}.
\end{proof}

Consider now a  weight $\P$-module $V$ such that $\N$ (and hence
$\N^{ps}$) acts trivially on~$V$ and~$c$ acts by a multiplication
by a nonzero scalar.  Then
\[
M_P(V)={\rm ind}(\P,\G;V)\simeq {\rm ind}\big(\m_{P}^{ps}\oplus \N_P^{ps}, \G; U(\m_P^{ps})V\big)
\]
and   we obtain immediately

\begin{corollary}\label{cor-induced-irred-m-type Ib}
Let $P$ be a parabolic subset of $\D$ such that $\G_{P}^0$ is finite-dimensional and $\m_{P}$ is
infinite-dimensional. Consider an irreducible weight $\G_{P}^0$-module
$V$ which is annihilated by~$\G_{P}^+$ and on which the central element $c$ acts injectively.
Then for any submodule $U$ of $M_P(V)$
 there exists a submodule $V_U$ of $U(\m_P)V$ such that
\[
U\simeq {\rm ind}(\m_{P}\oplus\N, \G; V_U).
\] Moreover,  $M_P(V)$ is irreducible if and only
if $U(\m_P)V$ is
 is irreducible.
\end{corollary}

Corollary~\ref{cor-induced-irred-m-type Ib} is essentially
Theorem~8 of \cite{FKM1}. In particular it reduces the case of
parabolic subalgebras of type Ib to the case of parabolic
subalgebras of type II.

\begin{example}\qquad
\begin{itemize}\itemsep=0pt
\item Examples of irreducible dense modules with non-zero central
charge were constructed in~\cite{CP} as tensor products oh highest
and lowest weight modules. Applying functor of pseudo parabolic
induction
 to these modules one obtains new examples of irreducible modules
over $\G$ with inf\/inite-dimensional weight spaces.

\item Series of
irreducible cuspidal modules over the Heisenberg subalgebra with a
non-zero central charge were constructed in~\cite{BBF}. These
modules have inf\/inite-dimensional weight spaces. We can not apply
functor of pseudo parabolic induction
 to these modules since the action of the Heisenberg subalgebra is
 torsion free. On the other hand, in the case of $A_1^{(1)}$, the
 functor of parabolic induction applied to such modules gives
 again new irreducible modules.
\end{itemize}
\end{example}

Parabolic induction can be easily generalized to the non-weight
case as follows (cf.~\cite{FKM1}). Let $\P=\P_0\oplus \N$ be a
parabolic subalgebra of type II, $P$ corresponding parabolic
subset of $\D$, $P^0=P\cap -P$. Let $\H'$ be the linear span of
$[\G_{\alpha}, \G_{-\alpha}]$, $\alpha\in P^0$ and $\H_P$ a
complement of~$\H'$ in~$\H$ such that $[\P_0, \H_P]=0$.

Let $\Lambda$ be an arbitrary Abelian category of $\P_0$-modules (note that $\Lambda$ may have a dif\/ferent
Abelian structure than the category of modules over $\P_0$). Given $V\in \Lambda$
and $\lambda\in \H_P^*$ one makes $V$ into a $\P$-module with $h|_{V}=\lambda(h)Id$ for any $h\in \H_P$ and $\N V=0$.
 Then one can construct a $\G$-module $M_P(V, \lambda)$ by parabolic induction.
It follows  from the construction that $M_P(V, \lambda)$ is
$\H_P$-diagonalizable. Analogously to Theorem~\ref{induced-verma-irred} one can show

\begin{theorem}\label{induced-verma-irred-nonweight}
If the central element $c$ acts injectively on $V\in \Lambda$ then for any submodule
$U$ of $M_P(V, \lambda)$ there exists a submodule $V_U$ of $V$ such that
\[
U\simeq M_P(V_U, \lambda).
\]
In particular, $M_P(V, \lambda)$
is irreducible if and only $V$ is irreducible.
\end{theorem}

\section{Categories of induced modules}\label{section4}

Let $\G$ be an af\/f\/ine Lie algebra of rank $> 1$, $\P=\P_0\oplus
\N$  a {\em pseudo} parabolic subalgebra of $\G$ of type II, $P$
corresponding parabolic subset of $\Delta$, $W(\P_0)$ the category
of weight (with respect to $\H$) $\P_0$-modules $V$ with an
injective action of $c$.

Denote by ${\mathcal O}(\G, \P)$ the category of weight $\G$-modules $M$ such that the  action of the
central element $c$ on $M$ is injective  and $M$ contains a nonzero $\P$-primitive element. Modules
$M_P(V)$ and $L_P(V)$ are typical objects of ${\mathcal O}(\G, \P)$.

For $M\in {\mathcal O}(\G, \P)$ we denote by $M^{\N}$ the subspace
of $M$ consisting of $\P$-primitive elements, that is the subspace
of $\N$-invariants. Clearly, $M^{\N}$ is a $\P_0$-module and hence
$M^{\N}\in W(\P_0)$. Let $\tilde{\mathcal O}(\G, \P)$ be the full
subcategory of ${\mathcal O}(\G, \P)$ whose objects $M$ are
generated by $M^{\N}$. Again $M_P(V)$ and $L_P(V)$ are objects of
$\tilde{\mathcal O}(\G, \P)$.

Both categories
${\mathcal O}(\G, \P)$ and $\tilde{\mathcal O}(\G, \P)$ are closed under the operations of taking submo\-dules,
quotients and countable direct sums.

 The parabolic induction provides a functor
\[
I: \ W(\P_0)  \to  {\mathcal O}(\G, \P), \qquad V  \mapsto
M_P(V)={\rm ind}(\P,\G; V).
\]
Note that $M_P(N)\simeq M_P^{ps}(N)$.
 The canonical image of $V$ in  $M_P(V)$ is annihilated by $\N$
  and, hence,
$I(V)$ is generated by its $\P$-primitive elements. Thus $I(V)\in \tilde{\mathcal O}(\G, \P)$.

In the opposite direction we have a  well
def\/ined functor
\[
R: \ {\mathcal O}(\G, \P) \to W(\P_0), \qquad M \mapsto M^{\N}.
\]

Denote by $\tilde{R}$ the restriction of $R$ onto the category $\tilde{\mathcal O}(\G, \P)$.

 We need the following lemma.

\begin{lemma} \label{lemma-ext}\qquad
\begin{itemize}\itemsep=0pt
\item Let $V\in W(\P_0)$, $M$ a subquotient of $M_P(V)$. Then
$w(M^{\N})\subset w(V)$. \item $M_P(V\oplus V')\simeq M_P(V)\oplus
M_P(V')$. \item Let $M\in \tilde{\mathcal O}(\G, \P)$ be generated
by $\P$-primitive elements then $M$ is a direct sum of modules of
type $M_P(V)$.
\end{itemize}
\end{lemma}

\begin{proof} Without loss of generality we will assume that $M$ is a quotient of $M_P(V)$. Suppose
 $M = M_P(V)/M'$. Then $M^{\N}\simeq V/(M'\cap V)$ (identifying $V$ with $1\otimes V$)
by Theorem \ref{induced-verma-irred}, and the f\/irst statement
follows. Second  statement follows immediately from the
def\/inition. Let $M\in \tilde{\mathcal O}(\G, \P)$ and
$M^{\N}=M'\oplus M''$ is a direct sum of $\P_0$-modules. If $M$ is
generated by $M^{\N}$ then $M\simeq M_P(M^{\N})$ by
Theorem~\ref{induced-verma-irred}, and hence $M\simeq
M_P(M')\oplus M_P(M'')$.
\end{proof}

{\samepage
\begin{theorem} \label{theorem-equivalence}\qquad
\begin{itemize}\itemsep=0pt
\item
 The functor $I: W(\P_0)  \to  {\mathcal O}(\G, \P)$ is a left adjoint to the
functor $R: {\mathcal O}(\G, \P) \to W(\P_0)$, that is $R \circ I$
is naturally isomorphic to the identity functor on $W(\P_0)$.

\item The functors $\tilde{R}$ and $I$ are mutually inverse equivalences of $W(\P_0)$
and the subcategory $\tilde{\mathcal O}(\G, \P)$.
\end{itemize}
\end{theorem}}

\begin{proof}
Let $V\in W(\P_0)$. Clearly, $V$ is naturally embedded into
$M_P(V)^{\N}$. On the other hand, $w(M_P(V)^{\N})\subset w(V)$ by
Lemma~\ref{lemma-ext}. Since $U(\P_0)V\simeq V$ then
$M_P(V)^{\N}\simeq V$.   If $M=M_P(V)$ then $R(M)\simeq V$ and
$(I\circ R)(M_P(V))\simeq M_P(V)$ by
Theorem~\ref{induced-verma-irred}. If $M$ is an arbitrary object
in ${\mathcal O}(\G, \P)$ then $M$ is generated by $M^{\N}=R(M)$
and $M^{\N}=\oplus_i M_i$, where $M_i$ are $\P_0$-modules and
$w(M_i)\cap w(M_j)=\varnothing$ if $i\neq j$. Then $M\simeq \oplus_i
M_P(M_i)$ by Theorem~\ref{induced-verma-irred}. On the other hand,
$I(M^(\N))\simeq \oplus_i M_P(M_i)$ by Lemma~\ref{lemma-ext}.
Hence, $(I\circ R)(M)\simeq M$, implying the statement.
\end{proof}

Theorem~\ref{induced-irred} is an immediate consequence of
Theorem~\ref{induced-verma-irred} or
Theorem~\ref{theorem-equivalence}.

Let $F(P)$ be a partial f\/lag of the parabolic subset $P$ of $\D$:
\[
F(P)=\{ \{0\} = F_k\subset
\cdots \subset F_1 \subset F_0=W \},
\]
$P^0=F_k\cap \Delta$, $\P_0=\G_P^0$, $\Delta^{\im}\subset F_k$. Then $\P_0$ is a subalgebra of $\G$
generated by  $\G_{\alpha}$, $\alpha\in F_k\cap \D$ and $\H$.

Fix $s$, $k\leq s< n$. Then $\delta\in F_s$ and the Lie subalgebra
$\m_{P}^s\subset\G$, generated by $\H$ and the root spaces
corresponding to the roots from $F_s\cap \D$, is
inf\/inite-dimensional. Obviously, it can be extended to a parabolic
subalgebra of $\G$ of type II, where $\m_{P}^s$ is the Levi
subalgebra: $\P_s=\m_P^s\oplus \N_s$, $\N_s\subset \N$.  If $V$ is
a $\P$-module with $\N N=0$ then $V'=U(\m_P^s)V$ is a
$\P_s$-module, $\N_s V'=0$. If $V'$ is a $\P_s$-module with a
trivial action of $\N_s$ and injective action of the central
element $c$ then the structure of the induced module ${\rm
ind}(\P_s, \G; V')$ is completely determined by the structure of
$V'$ by Theorem~\ref{induced-verma-irred}. In particular, this
module is irreducible if and only if~$V'$ is irreducible.
Moreover, since
\[
M_P(V)\simeq {\rm ind}(\P_s, \G; V'),
\] we have the following
interesting observation.

\begin{corollary}
If $V\in W(\P_0)$, $s$ is such that $k\leq s< n$ and
$V'=U(\m_P^s)V$ then the submodule structure of the induced module
$M={\rm ind}(\P_s, \G; V')$ is determined by the submodule
structure of $V$ $($as in Theorem~{\rm \ref{induced-verma-irred}}$)$.  In
particular, $M$ is irreducible if and only if $V$ is irreducible.
\end{corollary}

Denote by $W(\m_P^s)$ the category of weight $\m_P^s$-modules with injective action of the central element $c$
and by
${\mathcal O}(\G, \P_s)$ the category of weight $\G$-modules $M$ such that the  action of  $c$ on $M$ is
injective  and $M$ contains a nonzero $\P_s$-primitive element. Modules $M_{P_s}(V)$ and $L_{P_s}(V)$ are the
objects of ${\mathcal O}(\G, \P_s)$, $V\in W(\P_s)$. For $M\in {\mathcal O}(\G, \P_s)$ we denote by $M^{\N_s}$
the subspace of $M$ consisting of $\P_s$-primitive elements. Denote by $\tilde{\mathcal O}(\G, \P_s)$  the full
subcategory of ${\mathcal O}(\G, \P_s)$  whose objects $M$ are generated by $M^{\N_s}$.

Then we have
the following functors
\begin{gather*}
I_s: \ W(\m_P^s)  \to  {\mathcal O}(\G, \P_s), \qquad
V  \mapsto  M_{P_s}(V)={\rm ind}(\P_s,\G; V),
\\
R_s: \  {\mathcal O}(\G, \P_s) \to W(\m_P^s), \qquad M \mapsto M^{\N_s}
\end{gather*}
and $\tilde{R}_s$ the restriction of $R_s$ onto  $\tilde{\mathcal O}(\G, \P_s)$.

\begin{theorem} \label{theorem-equivalence-s}
For any  $s$, $k\leq s< n$,
\begin{itemize}\itemsep=0pt
\item
 the functor $I_s$ is a left adjoint to $R_s$;

\item the functors $\tilde{R}_s$ and $I_s$ are mutually inverse equivalences of $W(\m_P^s)$
and the subcategory $\tilde{\mathcal O}(\G, \P_s)$.
\end{itemize}
\end{theorem}

Note that for $k\leq s<r$  we have
$$
\xymatrix{
 {W(\m^s_P)\ } \ar@<0,5ex>[rr]^-{I_s}\ar@{^{(}->}[d]
&&  \tilde{\mathcal O}(\G,\P_s) \ar@<0,5ex>[ll]^-{\tilde{R}_s} \ar@{^{(}->}[d]\\
  {W(\m^r_P)\ } \ar@<0,5ex>[rr]^-{I_r } && \tilde{\mathcal O}(\G,\P_r) \ar@<0,5ex>[ll]^-{\tilde{R}_r }
           }
$$

Hence, $I_s$ is just the restriction of $I_r$ onto $W(\m_P^s)$,  while $\tilde{R}_s$ is the restriction
of $\tilde{R}_r$ onto  $\tilde{\mathcal O}(\G, \P_s)$.

\begin{remark}
One can establish similar category equivalences for non-weight
modules (cf.\ Section~5 in~\cite{FKM1}).
\end{remark}

\subsection*{Acknowledgements}
The f\/irst author is supported in part by the CNPq grant
(301743/2007-0) and by the Fapesp grant (2005/60337-2). The second
author is supported by the Fapesp grant (2007/025861).

\pdfbookmark[1]{References}{ref}
\LastPageEnding


\begin{thebibliography}{99}

\footnotesize\itemsep=0pt

\bibitem{BBF} Bekkert V., Benkart G., Futorny V., Weyl algebra
modules,  Kac--Moody Lie Algebras and Related Topics, {\it Contemp. Math.} {\bf 343} (2004), 17--42.

\bibitem{CP} Chari V., Pressley A., New unitary
representations of loop groups, {\it Math. Ann.} {\bf  275} (1986), 87--104.

\bibitem{Ch}  Chari V., Integrable representations of af\/f\/ine Lie algebras,
{\it Invent. Math.} {\bf  85} (1986),  317--335.

\bibitem{C}
Cox B., Verma modules induced from nonstandard Borel subalgebras,
{\it Pacific J. Math.} {\bf 165} (1994), 269--294.

\bibitem{DFP} Dimitrov I., Futorny V., Penkov I., A reduction theorem
for highest weight modules over toroidal Lie algebras, {\it Comm. Math. Phys.} {\bf 250} (2004), 47--63.

\bibitem{DG} Dimitrov I., Grantcharov D., Private communication.

\bibitem{DMP} Dimitrov I., Mathieu O., Penkov I., On the structure of weight modules,
{\it Trans. Amer. Math. Soc.} {\bf  352} (2000), 2857--2869.

\bibitem{DP}
Dimitrov  I., Penkov I., Partially integrable highest weight modules,
{\it Transform. Groups} {\bf 3} (1998), 241--253.

\bibitem{Fe}
Fernando S., Lie algebra modules with f\/inite-dimensional
 weight spaces.~I, {\it Trans. Amer.  Math. Soc.} {\bf  322} (1990), 757--781.

\bibitem{F1}  Futorny V., The weight
representations of semisimple f\/inite-dimensional Lie algebras,
in Algebraic Structures
and Applications, Kiev University,  1988, 142--155.

\bibitem{F2}
 Futorny V., Representations of af\/f\/ine Lie algebras, {\it Queen's Papers in Pure and Applied Mathematics}, Vol.~106, Queen's University, Kingston, ON, 1997.

\bibitem{F3} Futorny V., Irreducible non-dense
$A_1^{(1)}$-modules, {\it Pacific J.   Math.} {\bf  172} (1996), 83--99.

\bibitem{F4}
Futorny V., The parabolic subsets of root systems and corresponding
representations of af\/f\/ine Lie algebras,
in Proceedings of the International Conference on Algebra  (Novosibirsk, 1989),
{\it Contemp. Math.} {\bf 131} (1992), part~2, 45--52.

\bibitem{F5}  Futorny V., Imaginary
Verma modules for af\/f\/ine Lie algebras, {\it Canad. Math.  Bull.}
{\bf 37} (1994), 213--218.

\bibitem{FKM1} Futorny V., K\"onig S., Mazorchuk V., Categories of induced modules
for Lie algebras with triangular decomposition, {\it Forum Math.} {\bf 13} (2001), 641--661.

\bibitem{FKM2} Futorny V., K\"onig  S., Mazorchuk V., Categories of induced modules
and projectively stratif\/ied algebras, {\it Algebr. Represent. Theory}
{\bf 5} (2002), 259--276.

\bibitem{FS} Futorny  V., Saif\/i H.,
Modules of Verma type and new irreducible representations for
af\/f\/ine Lie algebras, in Representations of Algebras (Ottawa, ON, 1992), {\it CMS Conf. Proc.}, Vol.~14, Amer. Math. Soc., Providence, RI, 1993, 185--191.

\bibitem{FT}
Futorny V., Tsylke A., Classif\/ication of irreducible nonzero level modules
with f\/inite-dimensional weight spaces for af\/f\/ine Lie algebras, {\it J. Algebra} {\bf  238} (2001), 426--441.

\bibitem{JK}
Jakobsen  H.P., Kac V.G., A new class of unitarizable highest weight representations of
inf\/inite-dimensional Lie algebras, in Nonlinear Equations in Classical and Quantum Field Theory (Meudon/Paris, 1983/1984), {\it Lecture Notes in Phys.}, Vol.~226, Springer, Berlin, 1985, 1--20.

\bibitem{K}  Kac V.G., Inf\/inite-dimensional Lie algebras, 3rd ed.,
Cambridge University Press, Cambridge, 1990.

\bibitem{M} Mathieu O., Classif\/ication of irreducible weight modules, {\it Ann. Inst. Fourier (Grenoble)} {\bf  50} (2000), 537--592.

\bibitem{S} Spirin S., ${\mathbb Z}^2$-graded modules with one-dimensional components over the Lie algebra $A_1^{(1)}$,
{\it Funktsional. Anal. i Prilozhen.} {\bf 21} (1987), 84--85.

\end{thebibliography}
\end{document}